\documentclass[12pt]{amsart}
\usepackage{amscd,amsmath,amssymb,amsfonts}
\usepackage[cmtip, all]{xy}
\theoremstyle{plain}
\newtheorem{thm}{Theorem}

\theoremstyle{definition}

\numberwithin{thm}{section}
\numberwithin{equation}{section}

\newcommand{\ga}[2]{\begin{gather}\label{#1}#2 \end{gather}}

\newcommand{\surj}{\twoheadrightarrow}

\newcommand{\sX}{{\mathcal X}}

\newcommand{\A}{{\mathbb A}}

\newcommand{\C}{{\mathbb C}}

\newcommand{\F}{{\mathbb F}}

\renewcommand{\P}{{\mathbb P}}
\newcommand{\Q}{{\mathbb Q}}

\newcommand{\Z}{{\mathbb Z}}

\begin{document}

\title[cohomological divisibility]{
Cohomological divisibility and
point count divisibility}
\author{H\'el\`ene Esnault}
\address{
Universit\"at Essen, FB6, Mathematik, 45117 Essen, Germany}
\email{esnault@uni-essen.de}
\author{Nicholas M. Katz}
\address{Princeton University, Mathematics, Fine Hall, NJ 08544-1000, USA }
\email{nmk@math.princeton.edu}

\date{May 8, 2003}
\begin{abstract}
Let $X\subset \P^n$ be a closed scheme defined by $r$
homogeneous equations of degrees $d_1\ge d_2\ge \ldots \ge d_r$ over
the  finite field $\F_q$, with complement $U:=\P^n\setminus X$. Let
$\kappa$ be the maximum 
of $0$ and the integral part of the rational number $\frac{n-d_2-\ldots
-d_r}{d_1}$. We show that the eigenvalues of the geometric Frobenius acting
on the $\ell$-adic 
cohomology $H^i_c(U\times_{\F_q}\overline{\F_q}, \Q_\ell)$ with
compact supports are divisible by $q^\kappa$ as algebraic integers.

\end{abstract}
\maketitle
\begin{quote}

\end{quote}

\section{Introduction}
Nearly seventy years have passed since the theorem of
Chevalley-Warning  \cite{Ch}, \cite{War}: over a finite field $k=\F_q$ of
characteristic $p$,
 if $f\in k[X_1, \ldots, X_n]$
is a polynomial in $n\ge 1$ variables of degree $d<n$,
 then the number of solutions in $\A^n(k)$ of $f=0$
is divisible by $p$.
For $d\ge 1$, Ax \cite{Ax} improved this to divisibility by $q^\mu$, where
$\mu$ is the least non-negative integer which is $\ge \frac{n-d}{d}$. Thus
\ga{1.1}{\mu={\rm max}(0, \text{Ceiling}(\frac{n-d}{d})).}
In this form, the result remains true, if trivially so, without the assumption 
that $n > d$.
Ax also noted that, by an elementary inclusion-exclusion argument,
his result implied that if one had $r\ge 1$ polynomials $f_i\in
k[X_1, \ldots, X_n]$, with $f_i$ of degree $d_i\ge 1$, then the number of 
solutions in $\A^n(k)$ of the system of equations $f_1=\ldots=f_r=0$ is 
divisible by $q^\lambda$, where $\lambda$ is the least non-negative 
integer which is $\ge \frac{n-\sum_i d_i}{\sum d_i}$, i.e., 
\ga{1.2}{\lambda={\rm max}(0, \text{Ceiling}(\frac{n-\sum_i d_i}{\sum d_i})).}
Katz \cite{Ka} improved this to divisibility by
$q^\mu$, where $\mu$ is the least non-negative 
integer which is $\ge \frac{n-\sum_i d_i}{{\rm max}_i d_i}$, i.e.,
\ga{1.3}{\mu:=\mu(n; d_1,\ldots, d_r):=
{\rm max}(0, \text{Ceiling}(\frac{n-\sum_i d_i}{{\rm max}_i d_i})).} 
An elementary proof
of this result was given by Wan \cite{W0}.
Adolphson and Sperber \cite{AS} gave a divisibility estimate for
additive character sums which, they showed, included this result. Wan \cite{W}
gave a divisibility estimate for multiplicative character sums which,
he showed, included this result as well.

del Angel noted in \cite{P} that if one reverses the 
inclusion-exclusion argument which Ax used to pass from a result
for one polynomial to a result for several polynomials, then the estimate
\eqref{1.3} for several polynomial gives the following result
for one polynomial which is reducible:
if one has $r\ge 1$ polynomials $f_i\in k[X_1, \ldots, X_n]$, with $f_i$
of degree $d_i \ge 1$, then the number of solutions in $\A^n(k)$
of the single equation $f_1\cdots f_r=0$ is divisible by $q^{\mu}$,
for $\mu$ = $q^{\mu(n;d_1,\ldots,d_r)}$ as above.
If one replaces $f_1\cdots f_r$ by any product $f_1^{a_1}\cdots f_r^{a_r}$
with all $a_i\ge 1$, 
 one has the same solutions in $\A^n(k)$, so the
same divisibility. Notice that if one omits some of the factors in
$f_1^{a_1}\cdots f_r^{a_r}$, i.e., if one
allows some $a_i$ to be zero, 
then the number of solutions of
$f_1^{a_1}\cdots f_r^{a_r} = 0$ in $\A^n(k)$ remains divisible by the 
same  $q^{\mu}$ [simply because if one removes some 
of the $f_i$'s from consideration, the numerator 
$n-\sum_i d_i$ of $\mu$ increases, and its denominator ${\rm max}_i d_i$
decreases].

With this result in hand, del Angel did inclusion-exclusion a la Ax and got
the following result. Suppose one has a list of $r\ge 1$ polynomials
$f_i\in k[X_1,\ldots, X_n]$, with $f_i$ of degree $d_i \ge 1$. 
Take any $R\ge 1$,
and any list of $R\ge 1$ polynomials $g_j \in k[X_1,\ldots, X_n]$,
such that $g_j$ is a product $g_j=\prod_i f_i^{a_{i,j}}$, with all $a_{i,j}\ge 0
$. Then
the equations $g_1=\ldots=g_R=0$ define a closed subscheme $X\subset 
\A^{n}$, whose number of $k$-points
is divisible by $q^{\mu(n; d_1,\ldots, d_r)}$.

There is an obvious projective version of these affine results. If
the $f_i$ are all homogeneous, then all the $g_j$ are homogeneous, and 
the same equations $g_1=\ldots=g_R=0$  also define a closed subscheme 
$X^{{\rm proj}}\subset \P^{n-1}$ whose affine cone is the $X$ above. 
The number of $\F_q$-rational points of
$U:=\P^{n-1}\setminus X^{{\rm proj}}$ is easily seen to be
divisible by $q^{\mu(n; d_1,\ldots, d_r)}$. 

        The common feature of these results is that one has a 
separated $k$-scheme $X/k$ of finite type, and an integer $\kappa \ge 1$, 
such that for each finite extension $\F_{q^{\nu}}$ of $\F_q$, 
$|X(\F_{q^{\nu}})|$ is divisible by $q^{\nu \cdot \kappa}$.
Some years ago, Deligne posed the problem of showing 
that in any such situation $(X/k, \kappa)$, the divisibilities 
of numbers of points result from corresponding divisibilities 
of all the eigenvalues of Frobenius on all the compact $\ell$-adic 
cohomology groups, $\ell$ any prime other than $p$, of the variety $X/k$ 
in question. As explained by Ax \cite{Ax}, such 
divisibilities of numbers of points imply 
(and indeed are equivalent to) the corresponding divisibility 
of all reciprocal zeros and poles of the zeta function. 
But for any given $\ell \neq p$, the $\ell$-adic cohomological expression 
of the zeta function may well have cancellation: the problem 
is that there might be pairs of Frobenius eigenvalues on 
different $\ell$-adic groups which cancel, but which both have 
insufficient divisibility. 

        In the case of a general $(X/k, \kappa)$, Deligne's problem 
remains open. In this note, we show it has an affirmative 
solution when $(X/k, \kappa)$ arises as an instance of the 
divisibility theorems explained above. 

{\it Acknowledgements:}  We warmly
thank Spencer Bloch for his interest in
and for his incisive comments on the work
leading up to this article.

\section{Statement of results: the affine case}

We fix a finite field $k=\F_q$ of
characteristic $p$ and a prime number
$\ell \neq p$. For a separated scheme $X/k$  of finite type, we
abbreviate
 \ga{1.5}{H^i(X):=H^i(X\times_k \overline{k}, \Q_\ell), \
H^i_c(X):=H^i(X\times_k \overline{k}, \Q_\ell).
}
We denote by $F=F_k$ the geometric Frobenius endomorphism of these cohomology 
groups. By the Grothendieck-Lefschetz trace formula \cite{Gr}
\ga{1.6}{|X(k)|=\sum_i (-1)^i {\rm Trace}(F;H^i_c(X)), \ \text{and}\\
|X(k_\nu)|=\sum_i (-1)^i {\rm Trace}(F^\nu;H^i_c(X))\notag
}
for $k_\nu =\F_{q^\nu}$.
Recall \cite{DeInt}, Corollaire 5.5.3, that all 
the eigenvalues of $F$
acting on $H^i_c(X)$ are algebraic integers, i.e. elements of
$\overline{\Q_\ell}$ which are integral over $\Z$. 
\begin{thm} \label{mainthm}
Suppose one has a list of $r\ge 1$ polynomials $f_i\in k[X_1,\ldots, X_n]$,
with $f_i$ of degree $d_i\ge 1$. Take any integer $R\ge 1$, and any list of
$R\ge 1$ polynomials $g_j \in k[X_1,\ldots, X_n]$, such that each $g_j$ is a
product
 $g_j=\prod_i f_i^{a_{i,j}}$, with all $a_{i,j}\ge 0$. 
Denote by $X/k$ the closed
subscheme of $\A^n(k)$ defined by $g_1=\ldots =g_R=0$. 
Then every eigenvalue of
 $F$ on
every compact cohomology group $H^i_c(X)$ is divisible as an algebraic
integer by $q^{\mu(n;d_1,\ldots, d_r)}$.
\end{thm}
{F}rom the excision sequence
\ga{1.7}{\ldots \to H^i_c(\A^n \setminus X)\to H^i_c(\A^n)\to H^{i}_c(
 X)\to \ldots,}
and the well-known cohomological structure of $\A^n$, namely
\ga{1.8}{H^i_c(\A^n)=0 \ \text{for} \ i\neq 2n, \ H^{2n}_c(\A^n)=\Q_\ell(-n),}
where $\Q_\ell(-n)$ is the 1 dimensional
vector space over $\Q_\ell$ on which $F$ acts as $q^n$, one sees that Theorem
\ref{mainthm} is equivalent to
\begin{thm} \label{eqthm}
Under the hypotheses of Theorem \ref{mainthm}, every eigenvalue of $F$
on every compact cohomology group
$H^i_c(\A^n\setminus X)$ is divisible as an algebraic integer by
 $q^{\mu(n;d_1,\ldots, d_r)}$.
\end{thm}
Theorem \ref{eqthm} iself results from the following slightly sharper
result, suggested by \cite{DeInt}, Corollaire 5.3.3, (ii), and by 
\cite{W1}, \cite{EW}. For each integer $j \ge 0$, let us define
\ga{20}{\mu_j(n; d_1,\ldots, d_r):= j+\mu(n-j; d_1,\ldots, d_r).}
Thus
\ga{21}{\mu_0(n; d_1,\ldots, d_r) = \mu(n; d_1,\ldots, d_r),}
and we have the obvious inequalities
\ga{22}{\mu_{j+1}(n; d_1,\ldots, d_r) \ge \mu_j(n; d_1,\ldots, d_r).}

\begin{thm} \label{thmweight}
Under the hypotheses of Theorem \ref{mainthm}, we have the following results.
\begin{itemize}
\item[1)] Every eigenvalue of $F$ on every group
$H^i_c(\A^n\setminus X)$ is divisible as an algebraic integer by
$q^{\mu(n; d_1,\ldots, d_r)}$.
\item[2)] For each $j\ge 1$, every eigenvalue of $F$ on $H^{n+j}_c(\A^n\setminus
 X)$ is divisible as an algebraic integer by $q^{\mu_j(n;d_1,\ldots, d_r)}$.
\end{itemize}
\end{thm}

\section{The proof of Theorem \ref{thmweight}}
We first reduce to the case $R=1$. For $j$, with $1\le j\le R$, we denote by
$G_j \subset \A^n$ the hypersurface in $\A^n$ defined by $g_j=0$.
Then
\ga{2.1}{X=\cap _{j=1}^R G_j, \ \text{hence} \\
\A^n \setminus X
=\cup_{j=1}^R U_j, \ \text{with} \ U_j=\A^n\setminus G_j.\notag}
One considers the spectral sequence \cite{V}, \cite{DeCSP} 6.2.10.3,
\ga{2.2}{ E_1^{-a,b}=\oplus_{1\le j_1 <j_2<\ldots <j_{a+1}\le R}
H_c^b(\cap_{i=1}^{a+1} U_{j_i} ) \Rightarrow H^{b-a}_c(\A^n\setminus X),}
which is the cohomological incarnation of inclusion-exclusion.
Thus every eigenvalue of F on a given group 
$H^{b}_c(\A^n\setminus X)$ is an eigenvalue of $F$ on some group
$H_c^b(\cap_{i=1}^{a+1} U_{j_i} )$
for some $a\ge 0$. So it suffices to show the following two statements:
\begin{itemize}
\item[1)]For every $b \ge 0$ and every $a \ge 0$, every eigenvalue 
of $F$ on every group $H_c^b(\cap_{i=1}^{a+1} U_{j_i} )$
is divisible as an algebraic integer by
$q^{\mu(n;d_1,\ldots, d_r)}$
\item[2)] For every $j \ge 1$ and every $a \ge 0$, every eigenvalue 
of $F$ on every group $H_c^{n+j}(\cap_{i=1}^{a+1} U_{j_i} )$
is divisible as an algebraic integer by $q^{\mu_j(n;d_1,\ldots, d_r)}$.
\end{itemize}

But an $(a+1)$-fold intersection $\cap_{i=1}^{a+1} U_{j_i}$
is the complement of the hypersurface of equation
$\prod_{i=1}^{a+1}g_{j_i}=0$, which falls under the case $R=1$ of Theorem
\ref{thmweight}.
So it suffices to prove Theorem \ref{thmweight} universally for $R=1$. 
We use the weak Lefschetz theorem and induction on $n$, 
the number of variables. 
For $n=1$,  we have $\mu(1; d_1,\ldots, d_r) = 0$ 
(because $\sum_i d_i \ge 1$), so assertion 1) of the theorem is an 
(easy) instance of Deligne's integrality theorem \cite{DeInt}, Corollaire 5.5.3.
Assertion 2) of the theorem holds for $n = 1$ because the only 
possibly nonzero group $H^{1+j}_c(\A^1\setminus X)$ with $j \ge 1$ 
is $H^{2}_c(\A^1\setminus X)$, which is $\Q_\ell(-1)$. 
So suppose the theorem universally true for $n-1$, in the case $R=1$.

Whatever the nonempty hypersurface $G\subset \A^n$, $\A^n\setminus G$ is
smooth and geometrically connected of dimension $n$, so Poincar\'e
duality tells us that the cup product
\ga{2.3}{H^i_c(\A^n\setminus G)\times H^{2n-i}(\A^n\setminus G)\to H^{2n}_c(\A^n
\setminus G)\cong \Q_\ell(-n)}
is a perfect pairing of $F$-modules.
Since $\A^n\setminus G$ is affine of dimension $n$, the affine Lefschetz theorem
 tells us that
\ga{2.4}{H^i(\A^n\setminus G)=0 \ \text{for} \ i>n.}
So by duality we have
\ga{2.5}{H^i_c(\A^n\setminus G)=0 \ \text{for}\  i<n.}
We now take $G$ defined by $\prod_if_i^{a_i}=0$. 

We first prove assertion 1) of the theorem. As noted by Wan and 
exploited in \cite{BEL}, Introduction  and in \cite{E2}, given the divisibility
\ga{2.6}{|(\A^n \setminus G)(k_\nu)|=\sum_{i\ge n} (-1)^i {\rm Trace}(F^\nu;
H^i_c(\A^n\setminus G))}
by $q^{\nu\cdot \mu}$  for every $\nu\ge 1$, it suffices to prove that for all
but a single value of $i$ we have the asserted divisibility of
Frobenius eigenvalues on $H^i_c(\A^n\setminus G)$. We choose $i=n$.
So we need only show the divisibility by $q^{\mu(n;d_1,\ldots, d_r)}$
of Frobenius eigenvalues on $H^i_c(\A^n\setminus G)$ for $i>n$.

This divisibility is invariant under finite extension of the ground field.
For a sufficiently general affine linear hyperplane $A\subset \A^n$
defined by an equation of the form
\ga{2.7}{A: \ \sum_{i=1}^n \alpha_i X_i + \beta_i=0}
with coefficients in a  finite extension of $k$, (which, by extending $k$,
we may assume to be $k$ itself), the restriction maps on ordinary cohomology
\ga{2.8}{H^i(\A^n\setminus G)\to H^i(A\setminus A\cap G)}
are injective for $i<n$. This is an instance of the weak Lefschetz
theorem  (see \cite{KaACT}, 3.4.1, applied to $V=\A^n\setminus G$, with $\pi$
the inclusion and $f$ the function $0$.) The Poincar\'e dual of this injection
is a surjection
\ga{2.9}{ H^{i-2}_c(A\setminus A\cap G)(-1)\surj H^i_c(\A^n\setminus G)}
for $i>n$. For $A$ general, each $f_i|_A$ will be a polynomial of the same
degree $d_i$. So by induction on $n$, we know that every eigenvalue of $F$
on every group $H^{i-2}_c(A\setminus A\cap G)$ is divisible by $q^{\mu
(n-1;d_1,\ldots, d_r)}$. So every eigenvalue of $F$ on every group
$H^{i-2}_c(A\setminus A\cap G)(-1)$ is divisible by $
q^{\mu(n-1;d_1,\ldots, d_r)+1}$. But we have the inequality
\ga{2.10}{\mu (n-1;d_1,\ldots, d_r)+ 1 =\mu_1(n; d_1,\ldots, d_r)\ge \mu
(n;d_1,\ldots, d_r).}
This concludes the proof of assertion 1) of Theorem \ref{thmweight}.

This same induction proves assertion 2) of Theorem \ref{thmweight}.
 Indeed, for any integer $j\ge 1$, the induction 
shows that for $i=n+j$, every eigenvalue of $F$ 
on $H^{n+j}_c(\A^n\setminus G)$
is also an eigenvalue of $F$ on
\ga{2.11}{H^{(n-1)+(j-1)}(A\setminus G\cap A)(-1).}
By induction, every eigenvalue of $F$ on \eqref{2.11} is divisible by
$q^{1+\mu_{j-1}(n;d_1,\ldots, d_r)}$. And by definition we have
\ga{2.12}{1+\mu_{j-1}(n;d_1,\ldots, d_r)=\mu_{j}(n;d_1,\ldots, d_r).}
This concludes the proof of Theorem \ref{thmweight}.

\section{The projective case}

Here is the projective version. Suppose now that the 
$f_i$ are all homogeneous forms. Then each $g_j$ 
is again homogeneous, and the closed 
$X/k$ in $\A^n$ defined by $g_1=\ldots =g_R=0$
is the affine cone over a
closed $X^{{\rm proj}}$ in $\P^{n-1}$, 
defined by the same equations.

\begin{thm} \label{thmproj}
Under the hypotheses of Theorem \ref{mainthm}, suppose in addition
all the $f_i$ are homogeneous. Then we have the following results.
\begin{itemize}
\item[1)] Every eigenvalue of $F$ on every group
$H^i_c(\P^{n-1}\setminus X^{{\rm proj}})$ is 
divisible as an algebraic integer by
$q^{\mu(n; d_1,\ldots, d_r)}$.
\item[2)] For each $j\ge 1$, every eigenvalue of $F$ on 
$H^{n+j}_c(\P^{n-1}\setminus X^{{\rm proj}})$ is divisible as an 
algebraic integer by $q^{\mu_j(n;d_1,\ldots, d_r)}$.
\end{itemize}
\end{thm}
\begin{proof}  
To prove the theorem, repeat the proof of 
Theorem \ref{thmweight} 
essentially verbatim. First reduce to the 
case $R=1$ by the same spectral sequence 
argument. Then do induction on $n$. Again, 
the case $n=1$ holds trivially.
The open variety $\P^{n-1}\setminus X^{{\rm proj}}$ is affine, 
smooth, and geometrically connected of dimension $n-1$. 
The appropriate reference for the weak 
Lefschetz theorem in this context
is \cite{KaACT} Appendix (d'apr\`es Deligne), 
A5, applied with its 
$X :=\P^{n-1}\setminus X^{{\rm proj}}$, and $f$ 
the inclusion into $\P^{n-1}$.
\end{proof}

\section{Concluding remarks} \label{sec:rmks}

\subsection{ Concordance} 
In other work on this subject, e.g. in \cite{E}, \cite{ENS}, 
\cite{EW}, one sees a different but equivalent expression of 
the divisibility, in terms of the floor, or integral part [x] of
a real number x. One renumbers so that 
$d_1\ge d_2\ge \ldots \ge d_r$, and then one defines
\ga{5.1}{\kappa(n;d_1,\ldots, d_r) 
:= {\rm max}(0, [\frac{n-d_2-\ldots -d_r}{d_1}]).}
This $\kappa$ is related to $\mu$ by the identity
\ga{5.2}{\kappa(n;d_1,\ldots, d_r) = \mu(n+1;d_1,\ldots, d_r),}
whose verification we leave to the reader.

\subsection{ Axiomatization} The arguments we give here
could easily be axiomatized, and then applied 
to any pair of cohomology theories
$(H^\star, H_c^\star)$, defined on smooth 
separated $\overline{k}$-schemes
of finite type, with values in 
finite-dimensional vector spaces over a field $K$
of characteristic zero, with the following properties:
\begin{itemize}
\item[i)] zeta function formula for $H^\star_c$,
\item[ii)] Poincar\'e duality between $H^\star$ and  $H^\star_c$,
\item[iii)] affine Lefschetz and weak Lefschetz theorems for $H^\star$,
\item[iv)]  spectral sequence of an open covering for $H^\star_c$.
\end{itemize}
An initial run through our arguments with $\mu = 0$ 
shows that all Frobenius eigenvalues on $H^\star_c$ in any such theory
are algebraic integers. 
We expect that Berthelot's rigid cohomology will soon be such a pair of
theories. Indeed, (i) is proven in \cite{Et}, 
(ii) in \cite{BerCRAS}, and (iv) in 
\cite{BerInv}. Only (iii) is not yet proven. Once (iii) has been proven, 
then Theorems \ref{thmweight} and \ref{thmproj} 
will hold for the eigenvalues of Frobenius on compact 
rigid cohomology, with divisibility by $q^\mu$ 
as algebraic integers 
and hence with slopes of Frobenius $\ge \mu$.

\subsection{ Application to Hodge type} Deligne has also formulated a Hodge 
theoretic version of the problem considered here. 
Given a nonempty separated $\C$-scheme $X/\C$ of finite type,
recall that its compact support Hodge type $\alpha(X/\C)$ is the largest 
integer $\alpha$ such that the Hodge filtration ${\rm Fil}$, for the
mixed Hodge structure on $H^\star_{c}(X)$, satisfies
\ga{5.3}{{\rm Fil}^\alpha H^\star_{c}(X) = H^\star_{c}(X).}
We can also "spread out" $X/\C$ to a separated scheme $\sX/A$
of finite type, with $A$ a subring of $\C$ which is finitely generated
as a $\Z$-algebra. For each such spreading out $\sX/A$, we define
$\beta(\sX/A)$ to be the largest integer $\beta$ with the 
following property: for every finite field $\F_q$, and for every
ring homomorphism $\phi : A \rightarrow \F_q$,
the number of points on the scheme $\sX_{\phi}/\F_q$ obtained from $\sX/A$
by the extension of scalars $\phi$ is divisible by $q^\beta$.
We then define $\beta(X/\C)$ to be the maximum of $\beta(\sX/A)$
for all spreadings out of $X/\C$. 
One sees easily that $\beta(X/\C) \le \text{dim}(X)$. The problem is to prove
that $\alpha(X/\C) =\beta(X/\C)$. 
[Presumably this equality holds because $\alpha(X/\C)$
and $\beta(X/\C)$ are both equal to some intrinsic "motivic" 
invariant $\gamma(X/\C)$, cf. \cite{GrBr}, sections 9-10 
for Grothendieck's discussion of this sort of question.]
In general, not much is known. But in cases where the divisibility 
theorems of Ax et al give an estimate 
$\beta(X/\C) \ge \kappa$, one knows that $\alpha(X/\C) \ge \kappa$,
cf. \cite{DD},  \cite{ENS}, \cite{P}, \cite{EW}. 

Indeed, given 
the results on the Hodge type of hypersurfaces as input, 
the (axiomatizations of the) proofs we give here of Theorems \ref{thmweight}
and \ref{thmproj} also prove their analogues over $\C$, 
with eigenvalue divisibility replaced by Hodge type. Consider first
the projective case. The first assertion of 
the Hodge analogue of Theorem \ref{thmproj}
is the result \cite{DD}, \cite{E}, \cite{ENS}, \cite{P}. 
To prove the second assertion, reduce to
the $R=1$ case by the spectral sequence, and then use induction and
the weak Lefschetz theorem exactly as in the finite field case.
Consider now the affine case, i.e., the Hodge analogue of 
Theorem \ref{thmweight}. The spectral sequence reduces us to
the case $R=1$. Once we know the first assertion in the case $R=1$,
we get the second assertion in that case by using
induction and the affine weak Lefschetz
theorem, exactly as in the finite field case.
To prove the first assertion in the $ R=1$ case, we reduce it to
the projective case. For this reduction,
for each $f_i\in \C[X_1,\ldots, X_n]$,
with $f_i$ of degree $d_i\ge 1$, denote by 
$F_i\in \C[X_0,\ldots, X_n]$ the homogeneous form of the
same degree $d_i$ such that 
\ga{5.4}{f_i(X_1,\ldots, X_n) = F_i(1, X_1,\ldots, X_n),}
and denote
by $F_{i,0}\in \C[X_1,\ldots, X_n]$ the leading form of $f_i$,
i.e., 
\ga{5.5}{F_{i,0}(X_1,\ldots, X_n) = F_i(0, X_1,\ldots, X_n).} 
Given $G$ the hypersurface in $\A^n$ defined
by $\prod_if_i^{a_i}=0$, denote
by $G^{{\rm proj}}$ the hypersurface in $\P^n$ defined by
$\prod_iF_i^{a_i}=0$, and denote by $G_0^{{\rm proj}}$ the 
hypersurface in $\P^{n-1}$ defined by
$\prod_iF_{i,0}^{a_i}=0$. Then we have an excision sequence
\ga{5.6}{\ldots \to H^{i-1}_c(\P^{n-1} \setminus G_0^{{\rm proj}}) \to 
H^i_c(\A^n \setminus G)\to H^i_c(\P^n \setminus G^{{\rm proj}})\to 
\ldots.}
But $H^{\star}_c(\P^{n-1} \setminus G_0^{{\rm proj}})$ has Hodge type
at least $\mu(n;d_1,\ldots, d_r)$, and 
$H^{\star}_c(\P^{n} \setminus G^{{\rm proj}})$ has Hodge type
at least 
\ga{5.7}{\mu(n+1;d_1,\ldots, d_r) \ge \mu(n;d_1,\ldots, d_r).}
Thus $H^{\star}_c(\A^n \setminus G)$ has
Hodge type
at least $\mu(n;d_1,\ldots, d_r)$, as required.

\bibliographystyle{plain}

\renewcommand\refname{References}

\end{document}